\documentstyle{amsppt}
\magnification1095
\voffset-7.5mm
\hoffset-2.5mm
\pagewidth{38.33pc} \pageheight{52.5pc}
\baselineskip=12pt plus0.2pt minus0.2pt
\abovedisplayskip=7pt plus2pt minus2pt
\belowdisplayskip=7pt plus2pt minus2pt
\abovedisplayshortskip=5pt plus1pt minus3pt
\belowdisplayshortskip=5pt plus1pt minus3pt
\hfuzz=2.5pt\rightskip=0pt plus1pt
\binoppenalty=10000\relpenalty=10000\relax
\TagsOnRight
\loadbold
\nologo
\addto\tenpoint{\normalbaselineskip=1.2\normalbaselineskip\normalbaselines}
\addto\eightpoint{\normalbaselineskip=1.2\normalbaselineskip\normalbaselines}
\let\wt\widetilde
\let\le\leqslant
\let\ge\geqslant
\let\[\lfloor
\let\]\rfloor
\let\<\langle
\let\>\rangle
\let\epsilon\varepsilon
\let\phi\varphi

\redefine\d{\roman d}
\redefine\ord{\operatorname{ord}}

\define\ba{\boldkey a}
\define\bb{\boldkey b}
\define\bc{\boldkey c}
\define\fg{\frak g}
\define\fG{\frak G}
\define\fS{\frak S}
\topmatter
\title
Diophantine problems for $q$-zeta values
\footnotemark
\endtitle
\author
Wadim Zudilin \rm(Moscow)
\endauthor
\date
\hbox to100mm{\vbox{\hsize=100mm%
\centerline{E-print \tt math.NT/0206179}
\smallskip
\centerline{03 June 2002}
\centerline{{\sl Last revision\/}: 06 July 2002}
}}
\enddate
\address
\hbox to70mm{\vbox{\hsize=70mm%
\leftline{Moscow Lomonosov State University}
\leftline{Department of Mechanics and Mathematics}
\leftline{Vorobiovy Gory, GSP-2, 119992 Moscow, RUSSIA}
\leftline{{\it URL\/}: \tt http://wain.mi.ras.ru/index.html}
}}
\endaddress
\email
{\tt wadim\@ips.ras.ru}
\endemail
\endtopmatter
\footnotetext{An extension of the talk given at the conference
{\it Probl\`emes Diophantiens\/} (Luminy, CIRM, May 6--10, 2002).}
\rightheadtext{Diophantine problems for $q$-zeta values}
\leftheadtext{W.~Zudilin}
\document

\subhead
1. Introduction
\endsubhead
As usual, quantities depending on a number~$q$
and becoming classical objects as $q\to1$ (at least formally)
are regarded as $q$-{\it analogues\/} or $q$-{\it extensions}.
A possible way to $q$-extend the values of the Riemann zeta function
reads as follows (here $q\in\Bbb C$, $|q|<1$):
$$
\zeta_q(k)
=\sum_{n=0}^\infty\sigma_{k-1}(n)q^n
=\sum_{\nu=0}^\infty\frac{\nu^{k-1}q^\nu}{1-q^\nu}
=\sum_{\nu=0}^\infty\frac{q^\nu\rho_k(q^\nu)}{(1-q^\nu)^k},
\qquad k=1,2,\dots,
\tag1
$$
where $\sigma_{k-1}(n)=\sum_{d\mid n}d^{k-1}$ is the sum
of powers of the divisors and the polynomials $\rho_k(x)\in\Bbb Z[x]$
can be determined recursively by the formulae
$\rho_1=1$ and $\rho_{k+1}=(1+(k-1)x)\rho_k+x(1-x)\rho_k'$
for $k=1,2,\dots$ (see \cite{1, Part~8, Chapter~1, Section~8,
Problem~75} for the case $k=2$). Then the limit relations
$$
\lim\Sb q\to1\\|q|<1\endSb(1-q)^k\zeta_q(k)
=\rho_k(1)\cdot\zeta(k)
=(k-1)!\cdot\zeta(k),
\qquad k=2,3,\dots,
\tag2
$$
hold; the equality $\rho_k(1)=(k-1)!$ is proved
in~\cite{2, formula~(7)}.
The above defined $q$-zeta values \thetag{1} present several
new interesting problems in the theory of diophantine approximations
and transcendental numbers; these problems are extensions
of the corresponding problems for ordinary zeta values and
we state some of them in Section~3 of this note.
Our nearest aim is to demonstrate how some recent contributions
to the arithmetic study of the numbers $\zeta(k)$, $k=2,3,\dots$,
successfully work for $q$-zeta values. Namely, we mean
the hypergeometric construction of linear forms
(proposed in the works of E.\,M.~Nikishin~\cite{3},
L.\,A.~Gutnik~\cite{4}, Yu.\,V.~Nesterenko~\cite{5}) and
the arithmetic method (due to G.\,V.~Chudnovsky~\cite{6},
E.\,A.~Rukhadze~\cite{7}, M.~Hata~\cite{8})
accompanied with the group-structure scheme
(due to G.~Rhin and C.~Viola~\cite{9},~\cite{10}).
The next section contains new irrationality measures
of the numbers $\zeta_q(1)$ and $\zeta_q(2)$
for $q^{-1}=p\in\Bbb Z\setminus\{0,\pm1\}$,
and our starting point is the following table illustrating
a connection of some objects and their $q$-extensions
(here $\[\,\cdot\,\]$~denotes the integral part of a number
and the notation `l.c.m.' means the least common multiple).
We refer the reader to the book~\cite{11}
and the works~\cite{12}--\cite{14}, where a motivation
and a ground are presented.

\medskip
\line{\hss\vbox{\offinterlineskip
\halign to120.35mm{\strut\tabskip=1mm minus 1mm
\strut\vrule#&\hbox to50mm{\hss#\hss}&%
\vrule#&\hbox to70mm{\hss#\hss}&%
\vrule#\tabskip=0pt\cr\noalign{\hrule}
height9.2pt&ordinary objects&%
&$q$-extensions, $p=1/q\in\Bbb Z\setminus\{0,\pm1\}$&\cr
\noalign{\hrule\vskip1pt\hrule}
height15.2pt&numbers $n\in\Bbb Z$&%
&`numbers' $[n]_p=\dfrac{p^n-1}{p-1\vphantom{\big|}}\in\Bbb Z[p]$&\cr
\noalign{\hrule}
&primes $l\in\{2,3,5,7,\dots\}\in\Bbb Z$&%
&$\gathered
\text{irreducible reciprocal polynomials} \vphantom{\big|} \\ \vspace{-4.5pt}
\Phi_l(p)=\prod\Sb k=1\\(k,l)=1\endSb^l(p-e^{2\pi ik/l})
\in\Bbb Z[p]
\endgathered$&\cr
\noalign{\hrule}
&Euler's gamma function $\Gamma(t)$&%
&$\gathered
\text{Jackson's $q$-gamma function} \vphantom{\big|} \\ \vspace{-1.5pt}
\Gamma_q(t)=\frac{\prod_{\nu=1}^\infty(1-q^\nu)}
{\prod_{\nu=1}^\infty(1-q^{t+\nu-1})\vphantom{\big|_{0_0}}}\,(1-q)^{1-t}
\endgathered$&\cr
\noalign{\hrule}
&$\gathered
\text{the factorial $n!=\Gamma(n+1)$} \vphantom{\big|} \\ \vspace{-4pt}
n!=\prod_{\nu=1}^n\nu\in\Bbb Z
\endgathered$&%
&$\gathered
\text{the $q$-factorial $[n]_q!=\Gamma_q(n+1)$} \vphantom{\big|^0}
\\ \vspace{-3pt}
[n]_p!=\prod_{\nu=1\topsmash{\vphantom{|}}}^n\frac{p^\nu-1}{p-1}
=p^{n(n-1)/2}[n]_q!\in\Bbb Z[p]
\endgathered$&\cr
\noalign{\hrule}
height17pt&$\displaystyle
\ord_ln!=\biggl\[\frac nl\biggr\]
+\biggl\[\frac n{l^2\vphantom{|_{\displaystyle0_0}}}\biggr\]+\dotsb$&%
&$\displaystyle
\ord_{\Phi_l(p)}[n]_p!=\biggl\[\frac nl\biggr\]$, \
$l=2,3,4,\dots$&\cr
\noalign{\hrule}
&$\aligned
D_n&=\text{l.c.m.}(1,\dots,n) \\ \vspace{-3pt}
&=\prod_{\text{primes $l\le n$}}l^{\[\log n/\log l\]}
\in\Bbb Z
\endaligned$&%
&$\aligned
D_n(p)&=\text{l.c.m.}([1]_p,\dots,[n]_p) \vphantom{\big|^0} \\ \vspace{-4.5pt}
&=\prod_{l=1\topsmash{\vphantom{|}}}^n\Phi_l(p)\in\Bbb Z[p]
\endaligned$&\cr
\noalign{\hrule}
&the prime number theorem&&Mertens' formula&\cr
&$\displaystyle
\lim_{n\to\infty}\frac{\log D_n}n=1
$&&$\displaystyle
\lim_{n\to\infty}\frac{\log|D_n(p)|}{n^2\log|p|}
=\frac3{\pi^2\vphantom{|_{\displaystyle0_0}}}
$&\cr
\noalign{\hrule}
}}\hss}

If $\psi(x)$ is the logarithmic derivative of Euler's gamma function
and $\{x\}=x-\[x\]$~is the fractional part of a number~$x$,
then, for each demi-interval $[u,v)\subset(0,1)$,
Mertens' formula yields the limit relation
$$
\lim_{n\to\infty}\frac1{n^2\log|p|}
\sum_{l:\{n/l\}\in[u,v)}\log|\Phi_l(p)|
=\frac3{\pi^2}\bigl(\psi'(u)-\psi'(v)\bigr)
=\frac3{\pi^2}\int_u^v\d\bigl(-\psi'(x)\bigr)
\tag3
$$
(see \cite{15, Lemma~1}),
which can be regarded as a $q$-extension of the formula
$$
\lim_{n\to\infty}\frac1n
\sum\Sb\text{primes }l>\sqrt{Cn}\\\{n/l\}\in[u,v)\endSb\log l
=\psi(v)-\psi(u)
=\int_u^v\d\psi(x)
$$
in the arithmetic method \cite{6}--\cite{10}.

\subhead
2. Rational approximations to $q$-zeta values
and basic transformations
\endsubhead
Let $a_0$, $a_1$, $a_2$, and $b$ be positive integers
satisfying the condition $a_1+a_2\le b$.
Then, Heine's series
$$
F(\ba,b)
=\frac{\Gamma_q(b-a_2)}{(1-q)\Gamma_q(a_1)}
\sum_{t=0}^\infty
\frac{\Gamma_q(t+a_1)\,\Gamma_q(t+a_2)}
{\Gamma_q(t+1)\,\Gamma_q(t+b)}
\,q^{a_0t}
$$
becomes a $\Bbb Q(p)$-linear form $F(\ba,b)=A\zeta_q(1)-B$
with the property
$$
p^{-M}D_m(p)\cdot F(\ba,b)\in\Bbb Z[p]\zeta_q(1)+\Bbb Z[p];
\tag4
$$
here $M=M(\ba,b)$ is some (explicitly defined) integer and $m$~is
the maximum of the $6$-element set
$$
\gathered
c_{00}=a_0+a_1+a_2-b-1,
\qquad
c_{01}=a_0-1, \qquad c_{11}=a_1-1, \qquad c_{21}=a_2-1,
\\
c_{12}=b-a_1-1, \qquad c_{22}=b-a_2-1.
\endgathered
$$
Taking $H(\bc)=F(\ba,b)$ and using the stability of
the quantity
$$
\frac{F(a_0,a_1,a_2,b)}
{\Gamma_q(a_0)\,\Gamma_q(a_2)\,\Gamma_q(b-a_2)}
=\frac{H(\bc)}{\Pi_q(\bc)},
\qquad\text{where}\quad
\Pi_q(\bc)=[c_{01}]_q!\,[c_{21}]_q!\,[c_{22}]_q!
=p^{-N(\bc)}\Pi_p(\bc),
$$
under the action of the transformations
$$
\alignat2
\tau&=(c_{22} \; c_{21} \; c_{01} \; c_{11} \; c_{12} \; c_{00})\:
&(a_0,a_1,a_2,b)&\mapsto(a_1,b-a_1,a_0,a_0+a_2),
\\
\sigma&=(c_{11} \; c_{21})(c_{12} \; c_{22})\:
&(a_0,a_1,a_2,b)&\mapsto(a_0,a_2,a_1,b)
\endalignat
$$
we arrive at the better than~\thetag{4} inclusions
$$
p^{-M}D_m(p)\Omega^{-1}(p)\cdot F(\ba,b)\in\Bbb Z[p]\zeta_q(1)+\Bbb Z[p]
\tag5
$$
with
$$
\Omega(p)=\prod_{l=1}^m\Phi_l^{\nu_l}(p),
\qquad
\nu_l=\max_{\fg\in\<\tau^2,\sigma\>}\ord_{\Phi_l(p)}
\frac{\Pi_p(\bc)}{\Pi_p(\fg\bc)}.
\tag6
$$
In addition, trivial estimates for $F(\ba,b)$ and
explicit formulae for the coefficient~$A$ imply that
$$
|F(\ba,b)|=|p|^{O(b)},
\qquad
|A|\le|p|^{(a_0+a_1+a_2)b-(a_1^2+a_2^2+b^2)/2+O(b)}
\tag7
$$
with some absolute constant in~$O(b)$.

Note that the non-trivial transformation~$\tau$
of the quantity $H(\bc)/\Pi_q(\bc)$
has been obtained (in other notation) by E.~Heine
still in~1847. The transformation group
$\fG=\<\tau,\sigma\>$ of order~$12$ has no
ordinary analogue since corresponding (in limit $q\to1$)
Gau\ss's hypergeometric series are divergent.
We use the group $\<\tau^2,\sigma\>$
of order~$6$ instead of the total available group~$\fG$
to ensure the required condition $a_1+a_2\le b$. Now, choosing
$a_0=a_2=8n+1$, $a_1=6n+1$, and $b=15n+1$,
and taking in mind \thetag{5}, \thetag{7}, and \thetag{3}
we derive the following result.

\proclaim{Theorem 1}
For each $q=1/p$, $p\in\Bbb Z\setminus\{0,\pm1\}$,
the number~$\zeta_q(1)$ is irrational and its
irrationality exponent satisfies the estimate
$$
\mu(\zeta_q(1))\le2.42343562\dotsc.
\tag8
$$
\endproclaim

A value $\mu=\mu(\alpha)$ is said to be the {\it irrationality
exponent\/} of a real irrational number~$\alpha$ if $\mu$~is the least
possible exponent such that for any $\epsilon>0$ the inequality
$|\alpha-a/b|\le b^{-(\mu+\epsilon)}$
has only finitely many solutions in integers $a$ and~$b$.
The estimate~\thetag{8} can be compared with the previous result
$\mu(\zeta_q(1))\le2\pi^2/(\pi^2-2)=\allowmathbreak2.50828476\dots$
of P.~Bundschuh and K.~V\"a\"an\"anen in~\cite{12}
corresponding to the choice $a_0=a_1=a_2=n+1$ and $b=2n+2$
in the above notation.

Similar arguments with a simpler group $\<\sigma\>$ of order~$2$
can be put forward to improve W.~Van Assche's estimate
$\mu(\log_q(2))\le3.36295386\dots$ in~\cite{13}
for the following $q$-extension of~$\log(2)$:
$$
\log_q(2)=\sum_{\nu=1}^\infty\frac{(-1)^{\nu-1}q^\nu}{1-q^\nu}
=\sum_{\nu=1}^\infty\frac{q^\nu}{1+q^\nu}.
$$
Namely, in~\cite{14} we obtain the inequality
$\mu(\log_q(2))\le3.29727451\dots$
for $q^{-1}=p\in\Bbb Z\setminus\{0,\pm1\}$.

In the case of the numbers $\zeta_q(2)$, consider the
positive integers $(\ba,\bb)=(a_1,a_2,a_3,b_2,b_3)$
satisfying the conditions $a_j<b_k$, $a_1+a_2+a_3<b_2+b_3$
and the $q$-basic hypergeometric series
$$
\align
\wt F(\ba,\bb)
&=\frac{\Gamma_q(b_2-a_2)\,\Gamma_q(b_3-a_3)}{(1-q)^2\Gamma_q(a_1)}
\sum_{t=0}^\infty
\frac{\Gamma_q(t+a_1)\,\Gamma_q(t+a_2)\,\Gamma_q(t+a_3)}
{\Gamma_q(t+1)\,\Gamma_q(t+b_2)\,\Gamma_q(t+b_3)}
\,q^{(b_2+b_3-a_1-a_2-a_3)t}
\\
&=\wt A\zeta_q(2)-\wt B.
\endalign
$$
Then
$p^{-M}D_{m_1}(p)D_{m_2}(p)\cdot\wt F(\ba,\bb)
\in\Bbb Z[p]\zeta_q(2)+\Bbb Z[p]$,
where $m_1\ge m_2$ are the two successive maxima of
the $10$-element set
$$
c_{00}=(b_2+b_3)-(a_1+a_2+a_3)-1,
\qquad
c_{jk}=\cases
a_j-1 &\text{if $k=1$}, \\
b_k-a_j-1 &\text{if $k=2,3$},
\endcases
\quad j=1,2,3,
$$
and, in addition,
$$
|\wt F(\ba,\bb)|=|p|^{O(\max\{b_2,b_3\})},
\qquad
|\wt A|\le|p|^{b_2b_3-(a_1^2+a_2^2+a_3^2)/2+O(\max\{b_2,b_3\})}.
$$
The $\bc$-permutation group $\fG\subset\fS_{10}$ generated
by all permutations of $a_1,a_2,a_3$, the permutation of $b_2,b_3$,
and the permutation
$(c_{00} \; c_{22})(c_{11} \; c_{33})(c_{13} \; c_{31})$
has order~$120$ and is known in connection with the Rhin--Viola
proof~\cite{9} of the new irrationality measure for~$\zeta(2)$
(see also \cite{16, Section~6}).
In notation $\wt H(\bc)=\wt F(\ba,\bb)$, the quantity
$$
\frac{\wt H(\bc)}
{[c_{00}]_q!\,[c_{21}]_q!\,[c_{22}]_q!\,[c_{33}]_q!\,[c_{31}]_q!}
$$
is stable under the action of the group~$\fG$. This $\fG$-stability
yields the inclusions
$$
p^{-M}D_{m_1}(p)D_{m_2}(p)\wt\Omega^{-1}(p)\cdot\wt F(\ba,\bb)
\in\Bbb Z[p]\zeta_q(2)+\Bbb Z[p]
$$
with a quantity $\wt\Omega(p)$ defined like in~\thetag{6}.
Finally, choosing $a_1=5n+1$, $a_2=6n+1$, $a_3=7n+1$,
and $b_2=14n+2$, $b_3=15n+2$ we deduce the following result~\cite{17}.

\proclaim{Theorem 2}
For each $q=1/p$, $p\in\Bbb Z\setminus\{0,\pm1\}$,
the number~$\zeta_q(2)$ is irrational and its
irrationality exponent satisfies the estimate
$$
\mu(\zeta_q(2))\le4.07869374\dotsc.
\tag9
$$
\endproclaim

The quantitative estimates of type~\thetag{9} for $\zeta_q(2)$
have been not known before,
although the transcendence of $\zeta_q(2)$ for any algebraic
number~$q$ with $0<|q|<1$ follows from Nesterenko's theorem~\cite{18}.

It is nice to mention that the simpler choice of the parameters
$a_1=a_2=a_3=n+1$,
$b_2=b_3=2n+2$ also proves the irrationality of~$\zeta_q(2)$
for $q^{-1}\in\Bbb Z\setminus\{0,\pm1\}$, and the limit $q\to1$
produces Ap\'ery's original sequence~\cite{19} of rational approximations
to~$\zeta(2)$.

We would like to stress that using, like in~\cite{7}--\cite{10},
(multiple) $q$-integrals for the both series
$F(\ba,b)$ and $\wt F(\ba,\bb)$ in study of arithmetic
properties of the numbers $\zeta_q(1)$ and $\zeta_q(2)$
is in great difficulties. 
The reason of this is due to non-existance of a concept
of changing the variable of $q$-integration
(see \cite{20} and \cite{21, Section~2.2.4}).

\subhead
3. General problems for $q$-zeta values
\endsubhead
We start with mentioning that, for an even integer $k\ge2$,
the series $E_k(q)=1-2k\zeta_q(k)/B_k$, where $B_k\in\Bbb Q$
are Bernoulli numbers, is known to be the Eisenstein series.
Therefore the modular origin (with respect to the parameter
$\tau=\frac{\log q}{2\pi i}$) of the functions $E_4,E_6,E_8,\dots$
gives the algebraic independence of the functions
$\zeta_q(2),\zeta_q(4),\zeta_q(6)$ over $\Bbb Q[q]$,
while all other even $q$-zeta values are polynomials
in $\zeta_q(4)$ and $\zeta_q(6)$. In this sence, the
consequence of Nesterenko's theorem~\cite{18} ``the numbers
$\zeta_q(2),\zeta_q(4),\zeta_q(6)$ are algebraically
independent over~$\Bbb Q$ for algebraic $q$, $0<|q|<1$''
reads as a complete $q$-extension of the consequence
of Lindemann's theorem~\cite{22} ``$\zeta(2)$ is transcendental''.
Moreover, the transcendence of values of the function
$$
1+4\sum_{\nu=0}^\infty\frac{(-1)^\nu q^{2\nu+1}}{1-q^{2\nu+1}}
=\biggl(1+2\sum_{n=1}^\infty q^{n^2}\biggr)^2
\tag10
$$
at algebraic points~$q$, $0<|q|<1$, also follows from
Nesterenko's theorem (a proof of Jacobi's identity~\thetag{10}
can be found, e.g., in~\cite{23, Theorem~2});
the series on the left-hand-side of~\thetag{10} is
a $q$-analogue of the series
$$
4\sum_{\nu=0}^\infty\frac{(-1)^\nu}{2\nu+1}=\pi.
$$
The best known estimate for the irrationality exponent of~\thetag{10}
in the case $q^{-1}\in\Bbb Z\setminus\{0,\pm1\}$ is obtained
in~\cite{24}.

The limit relations~\thetag{2} as well as the expected algebraic
structure of the ordinary zeta values motivate the following
questions (we also regard $\zeta_q(1)$ to be an odd $q$-zeta value,
although the corresponding ordinary harmonic series is divergent).

\proclaim{Problem 1}
Prove that the $q$-zeta values $\zeta_q(1),\zeta_q(2),\zeta_q(3),\dots$
as functions of~$q$ are linearly independent over $\Bbb C(q)$.
\endproclaim

\proclaim{Problem 2}
Prove that the $q$-functional set involving the three even $q$-zeta values
$\zeta_q(2),\zeta_q(4),\zeta_q(6)$ and all
odd $q$-zeta values $\zeta_q(1),\zeta_q(3),\zeta_q(5),\dots$
consists of functions that
are algebraically independent over $\Bbb C(q)$.
\endproclaim

The associated diophantine problems consist in proving
the corresponding linear and algebraic independences over
the algebraic closure of~$\Bbb Q$ for algebraic~$q$
with $0<|q|<1$. In this direction,
even irrationality and $\Bbb Q$-linear independence
results for $q$-zeta values at the point $q\in\Bbb Q$
with $q^{-1}\in\Bbb Z\setminus\{0,\pm1\}$
would be very interesting.

A problem of other type is to construct a model of multiple
$q$-zeta values involving $q$-zeta values~\thetag{1}
and possessing similar properties with the model of
multiple zeta values~\cite{25}.



\Refs

\ref\no1
\by G.~P\'olya and G.~Szeg\"o
\book Problems and theorems in analysis
\vol2
\publaddr New York
\publ Springer-Verlag
\yr1976
\endref

\ref\no2
\by M.~Kaneko, N.~Kurokawa, and M.~Wakayama
\paper A variation of Euler's approach to values
of the Riemann zeta function
\jour Preprint (June 2002)
\finalinfo E-print {\tt math.QA/0206171}
\endref

\ref\no3
\by E.\,M.~Nikishin
\paper On irrationality of values of functions $F(x,s)$
\jour Mat. Sb. [Russian Acad. Sci. Sb. Math.]
\vol109
\yr1979
\issue3
\pages410--417
\endref

\ref\no4
\by L.\,A.~Gutnik
\paper On the irrationality of certain quantities involving~$\zeta(3)$
\jour Acta Arith.
\yr1983
\vol42
\issue3
\pages255--264
\endref

\ref\no5
\by Yu.\,V.~Nesterenko
\paper A few remarks on~$\zeta(3)$
\jour Mat. Zametki [Math. Notes]
\vol59
\yr1996
\issue6
\pages865--880
\endref

\ref\no6
\by G.\,V.~Chudnovsky
\paper On the method of Thue--Siegel
\jour Ann. of Math. (2)
\vol117
\issue2
\yr1983
\pages325--382
\endref

\ref\no7
\by E.\,A.~Rukhadze
\paper A lower bound for the approximation of~$\ln2$
by rational numbers
\jour Vestnik Moskov. Univ. Ser.~I Mat. Mekh.
[Moscow Univ. Math. Bull.]
\yr1987
\issue6
\pages25--29
\endref

\ref\no8
\by M.~Hata
\paper Legendre type polynomials and irrationality measures
\jour J. Reine Angew. Math.
\vol407
\issue1
\yr1990
\pages99--125
\endref

\ref\no9
\by G.~Rhin and C.~Viola
\paper On a permutation group related to~$\zeta(2)$
\jour Acta Arith.
\vol77
\issue1
\yr1996
\pages23--56
\endref

\ref\no10
\by G.~Rhin and C.~Viola
\paper The group structure for~$\zeta(3)$
\jour Acta Arith.
\vol97
\issue3
\yr2001
\pages269--293
\endref

\ref\no11
\by G.~Gasper and M.~Rahman
\book Basic hypergeometric series
\bookinfo Encyclopedia of Mathematics and its Applications
\vol35
\publaddr Cambridge
\publ Cambridge Univ. Press
\yr1990
\endref

\ref\no12
\by P.~Bundschuh and K.~V\"a\"an\"anen
\paper Arithmetical investigations of a certain infinite product
\jour Compositio Math.
\vol91
\yr1994
\pages175--199
\endref

\ref\no13
\by W.~Van Assche
\paper Little $q$-Legendre polynomials and irrationality
of certain Lambert series
\jour The Ramanujan J.
\yr2001
\vol5
\issue3
\pages295--310
\endref

\ref\no14
\by W.~Zudilin
\paper Remarks on irrationality of $q$-harmonic series
\jour Manuscripta Math.
\yr2002
\vol107
\issue4
\pages463--477
\endref

\ref\no15
\by E.~Heine
\paper Untersuchungen \"uber die Reihe $\dots$
\jour J. Reine Angew. Math. (Crelles J.)
\vol34
\yr1847
\pages285--328
\endref

\ref\no16
\by W.~Zudilin
\paper Arithmetic of linear forms involving odd zeta values
\jour Preprint (August 2001)
\finalinfo E-print {\tt math.NT/0206176}
\endref

\ref\no17
\by W.~Zudilin
\paper On the irrationality measure for $q$-analogue of~$\zeta(2)$
\jour Mat. Sb. [Russian Acad. Sci. Sb. Math.]
\vol193
\issue8
\yr2002
\endref

\ref\no18
\by Yu.\,V.~Nesterenko
\paper Modular functions and transcendence questions
\jour Mat. Sb. [Russian Acad. Sci. Sb. Math.]
\vol187
\issue9
\yr1996
\pages65--96
\endref

\ref\no19
\by R.~Ap\'ery
\paper Irrationalit\'e de~$\zeta(2)$ et~$\zeta(3)$
\jour Ast\'erisque
\vol61
\yr1979
\pages11--13
\endref

\ref\no20
\by R.~Askey
\paper The $q$-gamma and $q$-beta functions
\jour Appl. Anal.
\yr1978
\vol8
\pages125--141
\endref

\ref\no21
\by H.~Exton
\book $q$-Hypergeometric functions and applications
\bookinfo Ellis Horwood Ser. Math. Appl.
\publaddr Chichester
\publ Ellis Horwood Ltd.
\yr1983
\endref

\ref\no22
\by F.~Lindemann
\paper \"Uber die Zalh~$\pi$
\jour Math. Annalen
\vol20
\yr1882
\pages213--225
\endref

\ref\no23
\by G.\,E.~Andrews, R.~Lewis, and Z.-G.~Liu
\paper An identity relating a theta function to a sum of Lambert series
\jour Bull. London Math. Soc.
\vol33
\pages25--31
\yr2001
\endref

\ref\no24
\by T.~Matala-aho and K.~V\"a\"an\"anen
\paper On approximation measures of $q$-logarithms
\jour Bull. Austral. Math. Soc.
\vol58
\yr1998
\pages15--31
\endref

\ref\no25
\by M.~Waldschmidt
\paper Valeurs z\^eta multiples: une introduction
\jour J. Th\'eorie Nombres Bordeaux
\vol12
\yr2000
\pages581--595
\endref

\endRefs
\enddocument